\documentclass[11pt]{amsart}
 \usepackage{epsf}
\usepackage{epsfig}

 \title{Classification of Six-Point Metrics}

 \author{Bernd Sturmfels \ and  \ Josephine Yu }
 \address{Department of Mathematics, University of California,
              Berkeley, CA 94720}
 \email{[bernd,jyu]@math.berkeley.edu}

 \newtheorem{theorem}{Theorem}
 \newtheorem{lemma}[theorem]{Lemma}

 \newtheorem{remark}[theorem]{Remark}
 \newtheorem{corollary}[theorem]{Corollary}
 \newtheorem{proposition}[theorem]{Proposition}

 \theoremstyle{definition}
 
 \newtheorem{example}[theorem]{Example}

 \newcommand{\reals}{{\mathbb R}}

\begin{document}

\begin{abstract}
There are $339$ combinatorial types of
generic metrics on six points. They
correspond to the $339$ regular triangulations
of the second hypersimplex $\Delta(6,2)$, 
which  also has $14$ non-regular triangulations.
\end{abstract}


 \maketitle

 \section{The Metric Fan}

We consider the cone of all metrics
on the finite set $\, \{1,2,\ldots,n\}$:
$$ C_n \,\, = \,\, \quad
\bigl\{ \, d \in \reals^{\binom{n}{2}} \,\, : \,\,
d_{ij} \geq 0 \,\, \hbox{and} \,\,\,
d_{ij} + d_{jk} \geq d_{ik} ~~ \hbox{for all }~~ 1 \leq i,j,k \leq n
\,\bigr\}. $$
This is a closed convex pointed polyhedral cone. Its
extreme rays have been studied in 
combinatorial optimization \cite{DFMV, DL}.
Among the extreme rays are the \emph{splits}.
The splits are the metrics
$\,\sum_{i \in A} \sum_{j \not\in A} e_{ij} \,\in
\, \reals^{\binom{n}{2} } \,$
as $A$ ranges over nonempty subsets of $\{1,2,\ldots,n\}$.
There is an extensive body of knowledge
(see \cite{DL, Du}) also on the facets of the 
subcone of  $C_n$ generated by the splits.

Our object of study is a canonical subdivision of the metric cone $C_n$.
It is called the \emph{metric fan} and denoted ${MF}_n$.
 A quick way to define 
the metric fan ${MF}_n$ is to
say that it is the
\emph{secondary fan of the second hypersimplex}
$$ \Delta(n,2) \quad = \quad {\rm conv} 
\bigl\{ e_i + e_j \,\, : \,\, 1 \leq i < j \leq n \, \bigr\} \quad \subset \quad
\reals^n . $$
Every metric $d$ defines a regular polyhedral subdivision $\Delta_d$
of $\Delta(n,2)$ as follows. The vertices of $\Delta(n,2)$ are identified
with the edges of the complete graph $K_n$, and subpolytopes
of $\Delta(n,2)$ correspond to arbitrary subgraphs of $K_n$.
 A subgraph $G$ is a cell of $\Delta_d$ if there 
exists an $x \in \reals^n$ satisfying
\begin{equation*}
\hbox{$\, x_i + x_j = d_{ij} \,$ if $\{i,j \} \in G\,$ \ \ \ and 
\ \ \ $\, x_i + x_j > d_{ij} \,$ if $\{i,j\} \not\in G$.}
\end{equation*}
Two metrics $d$ and $d'$ lie in the same cone of the metric
fan $MF_n$ if they induce the same subdivision
$\,\Delta_d = \Delta_{d'}\,$ of the second hypersimplex $\Delta(n,2)$.
We say that the  metric $d$ is \emph{generic} if $d$ lies in an open
cone of $MF_n$. This is equivalent to saying that $\Delta_d$
is a \emph{regular triangulation} of $\Delta(n,2)$. 

These triangulation of $\Delta(n,2)$ and the resulting
metric fan $MF_n$ were studied by De Loera,
Sturmfels and Thomas \cite{DST}, who had been unaware of
an earlier appearance of the same objects
in phylogenetic combinatorics \cite{BD, Dre}. 
In \cite{Dre}, Dress considered the polyhedron
dual to the triangulation $\Delta_d$,
$$ P_d \quad = \quad \bigl\{
\,x \in {\reals}_{\geq 0}^n \,: \,\, x_i+x_j \,\geq \, d_{ij} \,\,\,
\hbox{for} \,\, 1 \leq i < j \leq n \, \bigr\}, $$
and he showed that its \emph{complex of bounded faces},
denoted $T_d$, is a natural object which generalizes
the phylogenetic trees derived from the metric $d$.
Both \cite{DST} and \cite{Dre} contain the description of
 the metric fans $MF_n$ for $n \leq 5$:
\begin{itemize}
\item The octahedron $\Delta(4,2)$ has three regular triangulations
$\Delta_d$. They are equivalent up to symmetry. The
corresponding \emph{tight span} $T_d$ is a quadrangle
with an edge attached to each of its four vertices. The 
three walls of the fan $MF_4$ correspond to the 
trees on $\{1,2,3,4\}$.
\item The fan $MF_5$ has $102$ maximal cones which
come in three symmetry classes. The tight spans $T_d$
of these three metrics are depicted in \cite[Figure A3]{Dre},
and the corresponding  triangulations $\Delta_d$ appear
(in reverse order) in \cite[page 414]{DST}. For instance, the
\emph{thrackle triangulation} of \cite[\S 2]{DST} corresponds to
the planar diagram in \cite{Dre}.
All three tight spans $T_d$
have five two-cells. (The type $T_{X,D_3}$ is slightly misdrawn in \cite{Dre}:
the two lower left quadrangles should form a flat pentagon).
\end{itemize}

The aim of this article
is to present the analogous classification for $n  = 6$.
The following result was obtained with the help of Rambau's software
TOPCOM \cite{Ram} for enumerating triangulations of arbitrary 
convex polytopes.

\begin{theorem} \label{main}
There are $194,160$ generic metrics on six points. These correspond
to the maximal cones in $MF_6$ and 
to the regular triangulations of $\Delta(6,2)$.
They come in $339$ symmetry classes. The hypersimplex
$\Delta(6,2)$ has also $3,840$ non-regular triangulations
which come in $14$ symmetry classes.
\end{theorem}

This paper is organized as follows. In Section 2
we describe all $12$ generic metrics whose tight span
$T_d$ is two-dimensional, and in Section 3 we describe
all $327$ generic metrics whose tight span has a
three-dimensional cell.  Similarly, in Section 4,
we describe the $14$ non-regular triangulations
of $\Delta(6,2)$. In each case a suitable system
of  combinatorial invariants will be introduced.
In Section 5 we study the geometry of the
metric fan $MF_6$.  The rays of $MF_6$ are precisely the
\emph{prime metrics} in \cite{KMT}. 
We determine the maximal cones incident to each prime metric,
and we discuss  the corresponding
minimal subdivisions of $\Delta(6,2)$. In Section 6 we present
a software tool for visualizing
the tight span  $T_d$ of any finite metric $d$.
This tool was written written
in POLYMAKE \cite{GJ}
with the help of Michael Joswig and Julian Pfeifle.
We also explain how its output  differs from 
the output of SPLITSTREE \cite{DHM}.

A complete list of all six-point metrics 
has been made available at
$$ \hbox{\tt bio.math.berkeley.edu/SixPointMetrics} $$
For each of the 339+14 types  in Theorem \ref{main},
the regular triangulation, Stanley-Reisner ideal,
and numerical invariants are listed.
The notation is consistent with that
used in the paper. In addition,
the webpage contains interactive pictures 
in JAVAVIEW \cite{java} 
of the tight span of each metric.

\medskip

\section{The $12$ Two-Dimensional Generic Metrics}

We identify each generic metric $d$ with its tight span $T_d$,
where the exterior segments have been contracted\footnote{Note that the exterior segments do appear in Figures 1--5 
of this paper and in the diagrams on our webpage.  They are drawn in green for extra clarity.} so that
every maximal cell has dimension $\geq 2$. With this
convention, generic four-point metrics are quadrangles
and five-point metrics are glued from five polygons
(cf.~\cite[Figure A3]{Dre}).
The generic six-point metrics, on the other hand, fall
naturally into two groups.

\begin{lemma}
Each generic metric on six points is either
a three-dimensional cell complex with
$26$ vertices, $42$ edges, $18$ polygons
and one $3$-cell, 
or it is a two-dimensional cell complex
with $25$ vertices, $39$ edges and $15$ polygons.
There are $327$ three-dimensional metrics
and $12$ two-dimensional metrics.
\end{lemma}

We first list the twelve types of two-dimensional metrics.
In each case the tight span consists of $15$ polygons
which are either triangles, quadrangles or pentagons.
Our first invariant is the vector $B = (b_3,b_4,b_5)$
where $b_i$ is the number of polygons with $i$ sides.
The next two invariants are the order of the symmetry
group and the number of cubic generators in the 
Stanley-Reisner ideal of the triangulation $\Delta_d$.
The last item is a representative metric
$ d  = (
d_{12},
d_{13},
d_{14},
d_{15},
d_{16},
d_{23},
d_{24},
d_{25},
d_{26},
d_{34},
d_{35},
d_{36},
d_{45},
d_{46},
d_{56}) $:

\smallskip

Type 1:
$(1, 10, 4), 1,2,
    (9, 9, 10, 13, 18, 18, 17, 6, 11, 17, 14, 9, 11, 8, 17)$

Type 2:
 $(1, 10, 4), 1,3, 
    (8, 8, 8, 14, 15, 16, 14, 6, 9, 12, 12, 7, 8, 7, 13)$

Type 3:
$ (1, 10, 4), 1,5, 
    (5, 6, 7, 8, 12, 11, 10, 5, 7, 11, 6, 6, 7, 5, 10)$

Type 4:
 $(1, 10, 4), 2,3, 
     (7, 5, 7, 12, 12, 12, 12, 5, 7, 10, 9, 7, 7, 5, 10)$

Type 5:
 $(1, 10, 4), 2,4, 
     (6, 7, 8, 10, 14, 13, 12, 6, 8, 13, 9, 7, 6, 6, 10)$

Type 6:
 $(1, 10, 4), 2,5, 
     (7, 7, 7, 11, 14, 12, 12, 6, 7, 14, 10, 7, 6, 7, 11)$

Type 7:
 $(1, 10, 4), 8,6, 
     (5, 5, 5, 8, 10, 10, 8, 5, 5, 8, 5, 5, 5, 5, 8)$
    
Type 8:
 $(2, 8, 5), 1,3 , 
     (5, 5, 7, 10, 11, 10, 10, 5, 8, 10, 7, 6, 5, 4, 7)$

Type 9:
 $(2, 8, 5), 2,4 , 
    (7, 7, 8, 10, 14, 14, 13, 5, 9, 13, 9, 7, 10, 6, 14)$

Type 10:
 $(2, 8, 5), 2,4 , 
     (5, 4, 5, 8, 9, 7, 8, 3, 6, 9, 6, 5, 5, 4, 7)$

Type 11:
 $(2, 8, 5), 2,4 , 
     (4, 5, 5, 8, 9, 9, 7, 4, 7, 8, 5, 4, 5, 4, 7)$

Type 12:
 $(3, 6, 6), 12,3, 
     (3, 3, 5, 6, 6, 6, 6, 3, 5, 6, 5, 3, 3, 3, 6)$

\smallskip

\noindent
The three metrics of types 9, 10 and 11 cannot be distinguished
by the given invariants. In Section 5 we explain how to 
distinguish these three types.

\begin{figure}[htb]
\epsfig{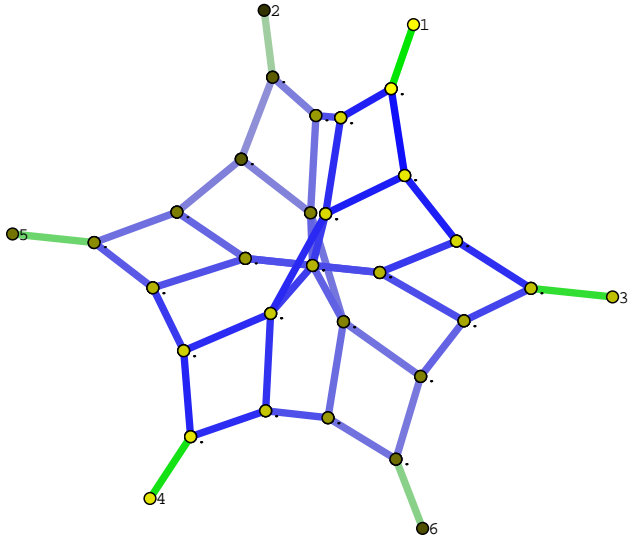}
\vskip -1.2cm
\caption{The tight span of the metric \# 12}
\label{fig:type12}
\end{figure}

The metric with the largest symmetry group is Type 12. Its
symmetry group has order $12$. This combinatorial type
of this metric is given by the Stanley-Reisner ideal
of the corresponding regular triangulation  of $\Delta(6,2)$:
 \begin{eqnarray*}  &\!\!\!\!\!\!\! \langle 
  x_{36} x_{14}, x_{25} x_{34}, x_{35} x_{46}, x_{16} x_{45},  x_{35} x_{12}, 
  x_{26} x_{35}, x_{36} x_{45}, x_{15} x_{36}, 
x_{26} x_{45}, x_{12} x_{46}, \\ &\!\!\!\!\!
  x_{12} x_{56}, x_{25} x_{36}, x_{45} x_{23}, x_{24} x_{13}, x_{45} x_{12}, 
 x_{34} x_{12}, x_{25} x_{46},  x_{23} x_{46}, x_{16} x_{25}, x_{13} x_{46}, \\ &\!\!\!\!\!
 x_{24} x_{36},  x_{35} x_{14}, x_{13} x_{56}, x_{26} x_{14}, x_{26} x_{13}, 
 x_{15} x_{46}, x_{36} x_{12}, x_{45} x_{13}, x_{25} x_{14}, x_{25} x_{13}, \\ &
 x_{15} x_{26} x_{34}, \,  x_{23} x_{56} x_{14},  \, x_{16} x_{24} x_{35} \,\rangle.
 \end{eqnarray*}
The number of quadratic generators is $30$, and this number is independent
of the choice of generic metric. The number of cubic generators of this
particular ideal is three (the last three generators), which is the 
third invariant listed under ``Type 12''. These cubic generators correspond
to ``empty triangles'' in the triangulation $\Delta_d$. For instance,
the cubic $x_{15} x_{26} x_{34}$ means that 
$\, {\rm conv} \{e_1+e_5, e_2+e_6, e_3+e_4\} \,$
is not a triangle in $\Delta_d$ but each of its three edges is
an edge in $\Delta_d$.
In the tight span $T_d$ this can be seen as follows:
\begin{eqnarray*}
 & \{ \hbox{geodesics between $1$ and $5$} \}\,\, \cap \,\,
 \{ \hbox{geodesics between $2$ and $6$} \} \\
  & \cap \quad
 \{ \hbox{geodesics between $3$ and $4$} \}
\quad  = \quad \emptyset ,
  \end{eqnarray*}
  but any two of these sets of geodesics have a common intersection.
  This can be seen in the picture of
  the tight span of the type 12 metric in Figure \ref{fig:type12}.
  
The twelve generic metrics listed above
demonstrate the subtle nature of the notion of
\emph{combinatorial dimension} introduced in  \cite{Dre}.
Namely, the combinatorial dimension of a generic metric $d$ 
can be less than that of a generic split-decomposable metric
\cite{BD}. This implies that the space of all 
$n$-point metrics of 
combinatorial dimension $\leq 2$ is a polyhedral fan
whose dimension exceeds the expected  number $4n-10 \,$
(cf.~\cite[Theorem 1.1 (d)]{DHM1}).

For six-point
metrics, this discrepancy can be understood by looking at the
centroid $\bigl(\frac{1}{3},
\frac{1}{3},\frac{1}{3},\frac{1}{3},\frac{1}{3},
\frac{1}{3}\bigr)$ of the hypersimplex $\Delta(6,2)$.
There are $25$ simplices in $\Delta(6,2)$
which  contain the centroid: the $15$ triangles
given by the perfect matchings of the graph $K_6$
and the $10$ five-dimensional simplices corresponding to two disjoint
triangles in $K_6$. In any given triangulation $\Delta_d$,
the centroid can lie in either one or the other. In the former
case, the tight span $T_d$ has a $3$-dimensional cell
dual to the perfect matching triangle in $\Delta_d$.
The combinatorial possibilities of these $3$-cells will
be explored in Section 3. In the latter case,
the tight span $T_d$ has a distinguished vertex
dual to the two-disjoint-triangles simplex in $\Delta_d$.
This vertex lies in nine polygons of $T_d$ which form
a link of type $K_{3,3}$. But there is no $3$-cell in $T_d$.
The distinguished vertex is the one in the center in Figure~1.

\medskip

\section{The $327$ Three-Dimensional Generic Metrics}

We next classify the $327$ three-dimensional metrics.
It turns out that in each case, the unique
$3$-dimensional cell is a simple polytope,
so its numbers $v$ of vertices and $e$ of edges
are determined by its number $f$ of faces:
$$ v = 2f - 4 \quad \hbox{and} \quad e = v+f-2 . $$
We consider the two vectors  $R = (r_3,r_4,r_5,r_6)$ and
$B = (b_3,b_4,b_5,b_6)$ where $r_i$ is the number
of polygons with $i$ edges on the $3$-cell
and $b_i$ is the number of polygons with $i$ edges not on the $3$-cell.
It turns out that no type has a polygon 
with $7$ or more sides. Hence the number
of facets of the $3$-cell is $\,f = r_3+r_4+r_5+r_6$.
Our third invariant is the pair $S = (s_2,s_3)$
where $s_2$ (resp.~$s_3$) is the number of $2/4$-splits
(resp.~$3/3$-splits) lying on the cone of the metric fan containing $d$.
The fourth invariant is the pair $C = (c_5,c_6)$ where
$c_i$ is the number of cubic generators of the
Stanley-Reisner ideal which involve $i$ of the points.
And finally we list $(g,t)$ where $g$ is the order
of the symmetry group and $t$ is the number of types
which share these invariants. These invariants
divides the $327$ three-dimensional
metrics into $251$ equivalence classes. We order the
classes lexicographically according to the vector
$(f,R,B,S,C,(g,t))$. The $251$ invariants 
are given in the following long list of strings
$\,\, R \,\, B \,\, S \,\, C \,\, gt \,$:
 \begin{smaller}
$$
\begin{matrix}
 & 4000 \, 0a22 \, 60 \, 02 \, 41
 & 4000 \, 0941 \, 50 \, 11 \, 11
 & 4000 \, 0941 \, 50 \, 21 \, 21
 & 4000 \, 1751 \, 40 \, 10 \, 11 \\
 & 4000 \, 1751 \, 40 \, 20 \, 11
 & 4000 \, 1832 \, 50 \, 11 \, 11
 & 4000 \, 1832 \, 50 \, 21 \, 11
 & 4000 \, 2642 \, 40 \, 20 \, 11 \\
 & 4000 \, 2642 \, 40 \, 20 \, 21
 & 4000 \, 2642 \, 40 \, 30 \, 21
 & 4000 \, 2642 \, 40 \, 40 \, 21
 & 4000 \, 2723 \, 50 \, 21 \, 21  \\
  & 4000 \, 3533 \, 40 \, 40 \, 11
 & 4000 \, 4424 \, 40 \, 60 \, 81
 & 2300 \, 0a21 \, 51 \, 21 \, 21
 & 2300 \, 0850 \, 40 \, 10 \, 11 \\
 & 2300 \, 0850 \, 50 \, 21 \, 11 
 & 2300 \, 0931 \, 40 \, 20 \, 11
 & 2300 \, 0931 \, 50 \, 11 \, 11
 & \underline{
   2300 \, 0940 \, 51 \, 01 \, 22}  \\
 & 2300 \, 0940 \, 51 \, 11 \, 13
 & 2300 \, 0940 \, 51 \, 21 \, 22
 & 2300 \, 0940 \, 61 \, 01 \, 21
 & 2300 \, 1660 \, 40 \, 20 \, 11 \\
 & 2300 \, 1660 \, 40 \, 40 \, 21
 & 2300 \, 1741 \, 40 \, 20 \, 12
 & 2300 \, 1741 \, 40 \, 30 \, 11
 & 2300 \, 1741 \, 50 \, 21 \, 11 \\ 
 & 2300 \, 1750 \, 51 \, 11 \, 13
 & 2300 \, 1750 \, 51 \, 21 \, 12
 & 2300 \, 1822 \, 40 \, 40 \, 21
 & 2300 \, 1831 \, 51 \, 21 \, 11 \\ 
 & 2300 \, 2551 \, 40 \, 40 \, 12
 & 2300 \, 2560 \, 51 \, 21 \, 22
 & 2300 \, 2632 \, 40 \, 40 \, 11
 & 2300 \, 2641 \, 51 \, 21 \, 21 \\
 & 2300 \, 3442 \, 40 \, 60 \, 21
 & 0600 \, 0a11 \, 41 \, 10 \, 11
 & 0600 \, 0a20 \, 42 \, 00 \, 22
 & 0600 \, 0a20 \, 52 \, 00 \, 12 \\
 & 0600 \, 0a20 \, 62 \, 00 \, 41
 & \underline{0600 \, 0c00 \, 63 \, 00 \, c1}
 & 0600 \, 0840 \, 41 \, 00 \, 11
 & 0600 \, 0840 \, 41 \, 00 \, 21 \\
 & 0600 \, 0840 \, 41 \, 10 \, 11
 & 0600 \, 0840 \, 41 \, 20 \, 21
 & 0600 \, 0840 \, 42 \, 00 \, 22
 & 0600 \, 0840 \, 51 \, 00 \, 21 \\
 & 0600 \, 0921 \, 41 \, 00 \, 21
 & 0600 \, 0930 \, 41 \, 00 \, 11
 & 0600 \, 0930 \, 41 \, 10 \, 11 
 & 0600 \, 0930 \, 51 \, 00 \, 12 \\
 & 0600 \, 0930 \, 51 \, 10 \, 11
 & 0600 \, 1650 \, 40 \, 30 \, 21
 & 0600 \, 1650 \, 40 \, 40 \, 21
 & 0600 \, 1650 \, 41 \, 10 \, 11 \\
 & 0600 \, 1650 \, 41 \, 20 \, 12
 & 0600 \, 1650 \, 42 \, 20 \, 11
 & 0600 \, 1731 \, 41 \, 20 \, 11
 & 0600 \, 1740 \, 41 \, 10 \, 11 \\
 & 0600 \, 1740 \, 41 \, 20 \, 11
 & 0600 \, 1740 \, 51 \, 10 \, 11
 & 0600 \, 1821 \, 41 \, 10 \, 11
 & 0600 \, 1821 \, 41 \, 30 \, 11 \\
 \end{matrix} $$ \end{smaller}  \begin{smaller}   $$ \begin{matrix}
 & 0600 \, 1830 \, 42 \, 20 \, 11
 & 0600 \, 1830 \, 52 \, 20 \, 11
 & 0600 \, 2460 \, 40 \, 40 \, 11
 & 0600 \, 2460 \, 41 \, 20 \, 21 \\
 & 0600 \, 2460 \, 42 \, 40 \, 41
 & 0600 \, 2541 \, 41 \, 40 \, 21
 & 0600 \, 2631 \, 41 \, 30 \, 11
 & 0600 \, 2640 \, 42 \, 40 \, 41 \\
 & 0600 \, 3270 \, 40 \, 60 \, 41
 & 2220 \, 0840 \, 40 \, 20 \, 11
 & 2220 \, 0840 \, 50 \, 21 \, 21
 & 2220 \, 0921 \, 40 \, 20 \, 21 \\ 
 & 2220 \, 0930 \, 51 \, 11 \, 13
 & 2220 \, 0930 \, 51 \, 21 \, 13
 & 2220 \, 1650 \, 40 \, 30 \, 21
 & 2220 \, 1650 \, 40 \, 40 \, 11 \\
 & 2220 \, 1650 \, 40 \, 40 \, 21
 & 2220 \, 1731 \, 40 \, 40 \, 11
 & 2220 \, 1740 \, 51 \, 21 \, 13
 & 2220 \, 2460 \, 40 \, 40 \, 11 \\
 & 2220 \, 2541 \, 40 \, 60 \, 21
 & 2220 \, 3270 \, 40 \, 60 \, 41
 & 0520 \, 0a01 \, 31 \, 00 \, 21
 & 0520 \, 0a10 \, 42 \, 00 \, 12 \\
 & 0520 \, 0a10 \, 52 \, 00 \, 11
 & 0520 \, 0b00 \, 53 \, 00 \, 21
 & 0520 \, 0830 \, 41 \, 00 \, 11
 & 0520 \, 0830 \, 41 \, 10 \, 11 \\
 & 0520 \, 0920 \, 31 \, 00 \, 11
 & 0520 \, 0920 \, 41 \, 00 \, 11
 & 0520 \, 0920 \, 41 \, 00 \, 21
 & 0520 \, 0920 \, 41 \, 10 \, 11 \\
 &  \underline{0520 \, 0920 \, 42 \, 00 \, 14}
 & 0520 \, 0920 \, 51 \, 00 \, 21
 & 0520 \, 0920 \, 52 \, 00 \, 21
 & 0520 \, 1640 \, 41 \, 10 \, 11 \\
 & 0520 \, 1640 \, 41 \, 20 \, 12
 & 0520 \, 1640 \, 42 \, 20 \, 11
 & 0520 \, 1730 \, 31 \, 20 \, 11
 & 0520 \, 1730 \, 41 \, 10 \, 11 \\
 & 0520 \, 1730 \, 41 \, 20 \, 11
 & 0520 \, 1730 \, 41 \, 30 \, 11
 & 0520 \, 1730 \, 42 \, 20 \, 12
 & 0520 \, 1811 \, 31 \, 20 \, 11 \\
 & 0520 \, 1820 \, 42 \, 20 \, 12
 & 0520 \, 1820 \, 52 \, 20 \, 11
 & 0520 \, 2450 \, 41 \, 40 \, 11
 & 0520 \, 2450 \, 42 \, 40 \, 21 \\
 & 0520 \, 2540 \, 41 \, 30 \, 11
 & 0520 \, 2621 \, 31 \, 40 \, 21
 & 0520 \, 2630 \, 42 \, 40 \, 21
 & 1330 \, 0830 \, 41 \, 10 \, 11 \\
 & 1330 \, 0830 \, 41 \, 20 \, 11
 & 1330 \, 0920 \, 41 \, 01 \, 11
 & 1330 \, 0920 \, 41 \, 10 \, 12
 & 1330 \, 0920 \, 41 \, 11 \, 11 \\
 & 1330 \, 0920 \, 51 \, 10 \, 11
 & 1330 \, 1640 \, 40 \, 40 \, 11
 & 1330 \, 1640 \, 41 \, 20 \, 13
 & 1330 \, 1640 \, 42 \, 20 \, 11 \\
 & 1330 \, 1730 \, 41 \, 11 \, 11
 & 1330 \, 1730 \, 41 \, 20 \, 11
 & 1330 \, 1730 \, 41 \, 30 \, 12
 & 1330 \, 1820 \, 42 \, 20 \, 11 \\
 & 1330 \, 1820 \, 52 \, 20 \, 11
 & 1330 \, 2450 \, 40 \, 60 \, 21
 & 1330 \, 2450 \, 41 \, 40 \, 11
 & 1330 \, 2450 \, 42 \, 40 \, 21 \\
 & 1330 \, 2540 \, 41 \, 30 \, 11
 & 1330 \, 2630 \, 42 \, 40 \, 21
 & 2221 \, 0920 \, 51 \, 21 \, 22
 & 2302 \, 0920 \, 51 \, 21 \, 21 \\
 & 3031 \, 1640 \, 40 \, 40 \, 11
 & 3031 \, 2450 \, 40 \, 60 \, 21
 & 0440 \, 0a00 \, 42 \, 00 \, 21
 & 0440 \, 0a00 \, 43 \, 00 \, 21 \\
 & 0440 \, 0910 \, 32 \, 00 \, 12
 & 0440 \, 0910 \, 41 \, 00 \, 11
 & 0440 \, 0910 \, 42 \, 00 \, 11
 & 0440 \, 1720 \, 31 \, 10 \, 12 \\
 & 0440 \, 1720 \, 31 \, 20 \, 14
 & 0440 \, 1720 \, 32 \, 10 \, 14
 & 0440 \, 1720 \, 32 \, 20 \, 11
 & 0440 \, 1720 \, 41 \, 10 \, 12 \\
 & 0440 \, 1720 \, 42 \, 20 \, 11
 & 0440 \, 1810 \, 32 \, 10 \, 14
 & 0440 \, 1810 \, 42 \, 10 \, 12
 & 0440 \, 1810 \, 42 \, 10 \, 22 \\
 & 0440 \, 1810 \, 42 \, 20 \, 11
 & 0440 \, 1900 \, 43 \, 10 \, 11
 & 0440 \, 2530 \, 31 \, 20 \, 12
 & 0440 \, 2530 \, 31 \, 40 \, 11 \\
 & 0440 \, 2530 \, 32 \, 30 \, 12
 & 0440 \, 2620 \, 32 \, 30 \, 12
 & 0602 \, 0a00 \, 42 \, 00 \, 21
 & 0602 \, 0a00 \, 43 \, 00 \, 41 \\
 & 0602 \, 0820 \, 42 \, 00 \, 41
 & 0602 \, 1720 \, 42 \, 20 \, 11
 & 0602 \, 2440 \, 42 \, 40 \, 41
 & 0602 \, 2620 \, 42 \, 40 \, 41 \\
 & 1331 \, 0820 \, 41 \, 10 \, 11
 & 1331 \, 0910 \, 41 \, 10 \, 11
 & 1331 \, 1720 \, 31 \, 20 \, 12
 & 1331 \, 1720 \, 41 \, 30 \, 11 \\
 & 1331 \, 1720 \, 42 \, 20 \, 12
 & 1331 \, 1810 \, 42 \, 20 \, 11
 & 1331 \, 2440 \, 41 \, 40 \, 11
 & 1331 \, 2440 \, 42 \, 40 \, 21 \\
 & 1331 \, 2530 \, 31 \, 40 \, 11
 & 1331 \, 2620 \, 42 \, 40 \, 21
 & 1412 \, 0910 \, 41 \, 11 \, 11
 & 2060 \, 2440 \, 42 \, 40 \, 41 \\
 & 2141 \, 0820 \, 41 \, 20 \, 21
 & 2141 \, 2620 \, 42 \, 40 \, 41
 & 2222 \, 1720 \, 41 \, 30 \, 11
 & 2222 \, 2440 \, 40 \, 60 \, 41 \\
 & 2222 \, 2440 \, 41 \, 40 \, 21
 & 4004 \, 2440 \, 40 \, 60 \, 81
 & 0360 \, 0900 \, 32 \, 00 \, 21
 & 0360 \, 0900 \, 33 \, 00 \, 61 \\
 & 0360 \, 1710 \, 31 \, 10 \, 11
 & 0360 \, 1710 \, 32 \, 20 \, 11
 & 0360 \, 2610 \, 22 \, 20 \, 13
 & 0360 \, 2610 \, 32 \, 20 \, 13 \\
 & 0360 \, 2700 \, 33 \, 20 \, 11
 & 0360 \, 2700 \, 33 \, 20 \, 21
 & 0360 \, 3420 \, 22 \, 40 \, 23
 & 0441 \, 0900 \, 31 \, 00 \, 21 \\
 & 0441 \, 1710 \, 31 \, 10 \, 11
 & 0441 \, 1710 \, 32 \, 10 \, 11
 & 0441 \, 1800 \, 32 \, 10 \, 12
 & 0441 \, 1800 \, 33 \, 10 \, 11 \\
 & 0441 \, 2520 \, 31 \, 20 \, 11
 & 0441 \, 2520 \, 32 \, 30 \, 11
 & 0441 \, 2610 \, 22 \, 20 \, 11
 & 0441 \, 2610 \, 32 \, 20 \, 11 \\
 & 0441 \, 2610 \, 32 \, 30 \, 11
 & 0441 \, 2700 \, 33 \, 20 \, 21
 & 0441 \, 3420 \, 22 \, 40 \, 21
 & 0522 \, 1710 \, 32 \, 10 \, 11 \\
 & 0522 \, 1800 \, 32 \, 10 \, 12
 & 0522 \, 1800 \, 33 \, 10 \, 11
 & 0522 \, 2520 \, 31 \, 20 \, 21
 & 0522 \, 2520 \, 32 \, 30 \, 11 \\
 & 0522 \, 2610 \, 32 \, 30 \, 11
 & 0603 \, 0900 \, 31 \, 01 \, 61
 & 1251 \, 1710 \, 31 \, 20 \, 11
 & 1251 \, 2520 \, 32 \, 30 \, 11 \\
 & 1332 \, 1710 \, 31 \, 20 \, 11
 & 1332 \, 1710 \, 32 \, 20 \, 11
 & 1332 \, 2520 \, 31 \, 40 \, 11
 & 1332 \, 2520 \, 32 \, 30 \, 11 \\
 & 1332 \, 2610 \, 32 \, 30 \, 12
 & 2304 \, 2520 \, 31 \, 40 \, 21
 & 0280 \, 2600 \, 22 \, 20 \, 11
 & 0280 \, 2600 \, 23 \, 20 \, 21 \\
 & 0280 \, 3410 \, 22 \, 30 \, 21
 & 0280 \, 3410 \, 22 \, 40 \, 21
 & 0361 \, 2600 \, 22 \, 20 \, 11
 & 0361 \, 2600 \, 23 \, 20 \, 11 \\
 & 0361 \, 3410 \, 22 \, 30 \, 11
 & 0361 \, 3410 \, 22 \, 40 \, 22
 & 0361 \, 3500 \, 23 \, 30 \, 11
 & 0442 \, 2600 \, 22 \, 20 \, 13 \\
 & 0442 \, 2600 \, 22 \, 20 \, 22
 & 0442 \, 2600 \, 23 \, 20 \, 22
 & 0442 \, 3410 \, 22 \, 30 \, 11
 & 0442 \, 3410 \, 22 \, 30 \, 23 \\
 & 0442 \, 3500 \, 23 \, 30 \, 13
 & 0523 \, 2600 \, 22 \, 20 \, 21
 & 0523 \, 2600 \, 23 \, 20 \, 21
 & 0523 \, 3410 \, 22 \, 40 \, 21 \\
 & 0604 \, 2600 \, 22 \, 20 \, 41
 & 1252 \, 3410 \, 22 \, 40 \, 22
 & 1333 \, 3410 \, 22 \, 40 \, 21
 & 1414 \, 3410 \, 22 \, 40 \, 21 \\
 & 0281 \, 4300 \, 13 \, 40 \, 12
 & 0281 \, 4300 \, 13 \, 40 \, 21
 & 0362 \, 4300 \, 13 \, 40 \, 11
 & 0362 \, 4300 \, 13 \, 40 \, 21 \\
 & 0443 \, 4300 \, 13 \, 40 \, 13
 & 0443 \, 4300 \, 13 \, 40 \, 21
 & 0524 \, 4300 \, 13 \, 40 \, 21
 & 0282 \, 6000 \, 04 \, 60 \, 41 \\
 & 0363 \, 6000 \, 04 \, 60 \, 21
 & 0363 \, 6000 \, 04 \, 60 \, 61
 & 0444 \, 6000 \, 04 \, 60 \, 83 \end{matrix}
 $$
 \end{smaller}

The letters $a,b$ and $c$ appearing in this list
represent the integers $10,11$ and $12$.
We label the $327$ types of three-dimensional 
metrics as Type 13, Type 14, \ldots, Type 339, 
in the order in which they appear in this list.
Whenever the string does not end in a $1$
then that string refers to more than one type.

\begin{figure}[htb]
\vskip -0.4cm
\epsfig{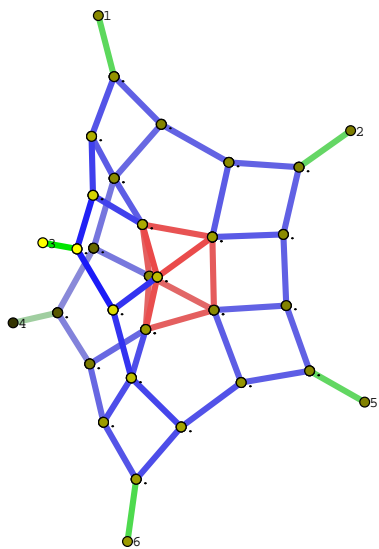}
\vskip -0.3cm
\epsfig{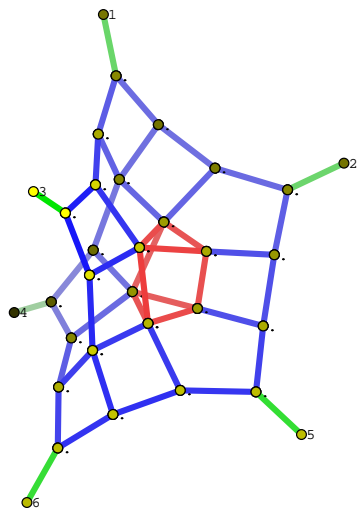}
\vskip -0.9cm
\caption{The tight spans of the metrics \# 32 and \# 33}
\label{fig:type32_33}
\end{figure}

For instance, the first underlined string
refers to Type 32 and Type 33. For 
both of these types, the three-dimensional
cell in $T_d$ is a triangular prism, hence
$\,R = (2,3,0,0)$. The next invariant
$\,B = (0,9,4,0)\,$ says that the two-dimensional
part of $T_d$ consists of nine quadrangles and 
four pentagons. Types 34 through 38 share
these characteristics. What distinguishes
Types 32-33 from Types 34-38 is the number of 
cubic generators in the Stanley-Reisner ideal.
The relevant vectors $C = (c_5,c_6)$ in the
table entries are $\,01\,$, $\,11 \,$ and $\,21$.
The tight spans of Type 32 and 33 are depicted in Figure \ref{fig:type32_33}.
The location of the four pentagons relative to the 
six exterior segments of the figure
 shows that these two types are non-isomorphic.

The second underlined string in the long list is 
Type 66. Its $3$-cell is a cube (hence $\,R = \,0600\,$), 
and its $2$-dimensional part consists of twelve quadrangles
(hence $\, B = \,0c00 \,$). This is  the unique generic metric 
which is \emph{split-decomposable} 
(hence $\, S = \, 63 \,$) in the sense of \cite{BD}.
It corresponds to the quadratic 
Gr\"obner bases  (hence $\, C = \,00 \,$)
and the \emph{thrackle triangulation} described in
\cite[\S 2]{DST}. It has the symmetry group
of a regular hexagon  (hence $g=12$) and 
it is uniquely characterized by $R$ and $B$
(hence $t = 1$). Its tight span is the logo
for the conference on \emph{Phylogenetic Combinatorics}
to be held in Uppsala, Sweden, in  July 2004,
{\tt http://www.lcb.uu.se/pca04/}.

\begin{figure}[htb]
\vskip -0.5cm
\epsfig{file = 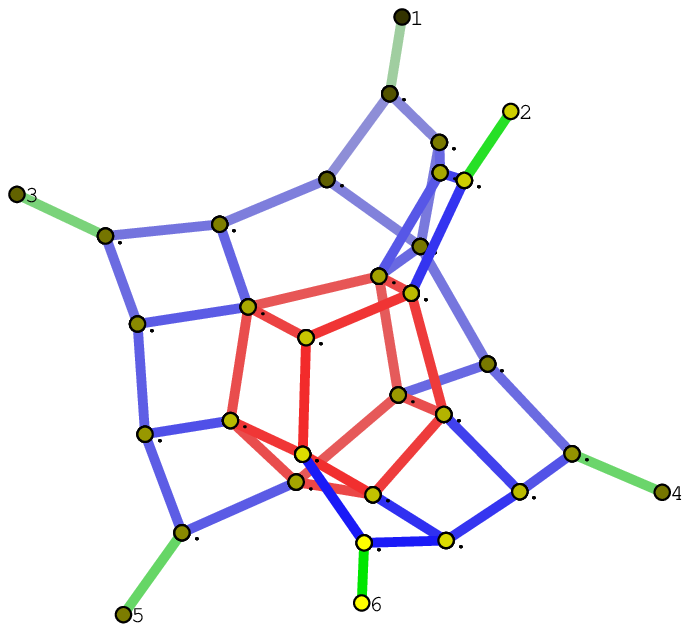, scale = 0.6}
\vskip -1.6cm
\caption{The tight span of the metric \# 131}
\label{fig:type131}
\end{figure}

The third underlined string represents a class of four types,
namely, Types 131, 132, 133 and 134.
In each of these four cases, the $3$-dimensional cell is
a pentagonal prism (hence $\, R = \,0520\,$),
the two-dimensional part consists of nine quadrangles
and two pentagons (hence $\, B = \,0920 \,$),
and the Gr\"obner basis is quadratic (hence $\, C = \,00 \,$).
Figure \ref{fig:type131} shows one of these tight spans.

\medskip

\section{The 14 Non-regular Triangulations}

Theorem 4.2 in \cite{DST} states that
the hypersimplex $\Delta(n,2)$ has non-regular
triangulations for $n \geq 9$. In this section we
strengthen this result as follows.

\begin{theorem}
The second hypersimplex $\Delta(n,2)$ admits
non-regular triangulations if and only if $n \geq 6$.
There are precisely $14$ symmetry classes of non-regular
triangulations of $\Delta(6,2)$.
\end{theorem}

It can be shown by explicit computations that all
triangulations of $\Delta(4,2)$ and $\Delta(5,2)$ are
regular. In what follows we list the $14$ non-regular
triangulations of $\Delta(6,2)$. Each of them can be
lifted to a non-regular triangulation of $\Delta(n,2)$ for 
$n \geq 7$ using the technique described at the end of \cite[\S 4]{DST}.

Each non-regular triangulation $\Delta$ has a dual 
polyhedral cell complex $T$. This complex $T$ shares all
the combinatorial properties of a tight span $T_d$, but
it cannot be realized as the complex of bounded faces
of a polyhedron $P_d$.  We call $T$ the
\emph{abstract tight span} dual to $\Delta$.

We use the  labels Type 340, Type 341,  $\ldots $, Type 353
to denote the $14$ non-regular triangulations $\Delta$ of $\Delta(6,2)$.
In each case we describe the abstract tight span $T$.
The first four of the $14$ abstract tight spans are
two-dimensional. 
They can be characterized by means of the invariants in Section 2:

\smallskip

Type 340: $\, (0, 12, 3), \,\,  1,\,  1 $

Type 341: $\, (0, 12, 3), \,\,  1,\, 2  $

Type 342: $\, (0, 12, 3), \,\,  6,\, 0  $

Type 343: $\, (1, 10, 4), \,\,  4,\, 1 $

\smallskip

\noindent The remaining ten abstract tight spans have a unique 
three-dimensional cell. We characterize them using the
invariants $(R,B,S,C,g)$ of Section 3:
 
\smallskip

Type 344: $ (4, 0, 0, 0), \, (0, 8, 6, 0), \, (4, 0),  \, (0, 0),\,   4$

Type 345: $ (0, 4, 4, 0), \,  (0, 8, 2, 0), \,(4, 0), \,  (0, 0), \,  4 $

Type 346: $ (0, 4, 4, 0), \,  (2, 4, 4, 0), \,(4, 0), \, (4, 0),   \, 2 $

Type 347: $ (0, 4, 4, 0), \,  (2, 4, 4, 0), \,(4, 0),\, (2, 0), \, 4 $

Type 348: $ (0, 4, 4, 0), \,  (2, 4, 4, 0), \, (4, 0),\, (2, 0),\, 4 $

Type 349: $ (0, 4, 4, 0), \,  (2, 4, 4, 0), \, (4, 0), \, (6, 0), \, 8 $

Type 350: $ (0, 3, 6, 0), \,  (2, 5, 2, 0), \, (3, 0), \, (2, 0), \, 2 $

Type 351: $ (0, 3, 6, 0), \,  (2, 5, 2, 0), \,  (3, 0), \, (4, 0), \, 2 $

Type 352: $ (0, 2, 8, 0), \,  (2, 6, 0, 0), \,  (2, 2), \, (2, 0), \,  4 $

Type 353: $ (0, 0, 12, 0), \, (6, 0, 0, 0),  \, (0, 4), \, (6, 0), \, 24 $.

\begin{figure}[htb]
\vskip -0.6cm
\epsfig{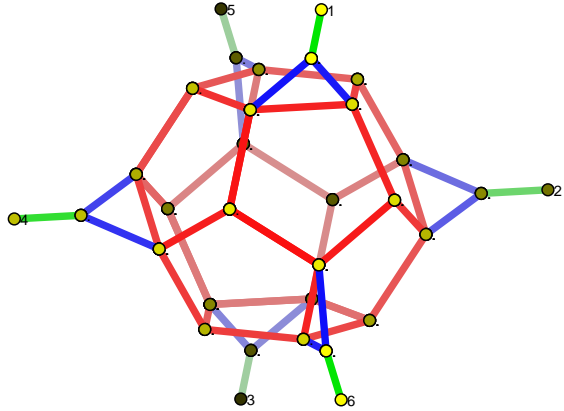}
\label{fig:dodec}
\vskip -1.8cm
\caption{The dodecahedral (abstract) tight span \# 353 }
\end{figure}

\noindent Type \# 353 is the most symmetric one
among triangulations $\Delta$ of $\Delta(6,2)$.
The corresponding abstract tight span $T$  is a 
beautiful object, namely, it is a dodecahedron with
six triangles and six edges attached, as shown in Figure \ref{fig:dodec}.
The Stanley-Reisner ideal corresponding to the dodecahedral type
\# 353 is generated by $30$ quadrics and $6$ cubics. The quadrics
in this ideal are
\begin{eqnarray*} &
x_{12} x_{35} , 
x_{12} x_{36} , 
x_{12} x_{45} , 
x_{12} x_{46} , 
x_{12} x_{56} , 
x_{13} x_{24} , 
x_{13} x_{26} , 
x_{13} x_{45} , \\ & 
x_{13} x_{46} , 
x_{13} x_{56} ,
x_{14} x_{23} , 
x_{14} x_{26} , 
x_{14} x_{35} , 
x_{14} x_{36} , 
x_{14} x_{56} , 
x_{15} x_{23} , \\ &
x_{15} x_{24} , 
x_{15} x_{26} , 
x_{15} x_{36} , 
x_{15} x_{46} , 
x_{23} x_{45} , 
x_{23} x_{46} , 
x_{23} x_{56} , 
x_{24} x_{35} , \\ &
x_{24} x_{36} , 
x_{24} x_{56} , 
x_{26} x_{35} , 
x_{26} x_{45} , 
x_{35} x_{46} , 
x_{36} x_{45} .
\end{eqnarray*}
The twelve variables $x_{ij}$ appearing
in this list can be identified with the
edges of an octahedron. The $30$ quadrics
are precisely the pairs of disjoint edges 
of the octahedron. We note that these quadratic 
generators (and hence the global structure of
the tight span) are also shared by the  last six
types (334, 335, 336, 337, 338, 339)
in the table of Section 3. The simplicial complex
represented by these $30$ quadrics is (essentially)
the boundary of the truncated octahedron
(with the six square faces regarded as tetrahedra).

Now, each of the Types 334, 335, 336, 337, 338, 339
and 353 has six cubic generators in its ideal. The choice
of these cubic generators amounts to subdividing each of
the six square faces of the truncated octahedron with
one of its two diagonals. Type 353 arises
from the most symmetric choice of these diagonals.
The six cubics in the ideal for Type 353 are
$$
\underline{x_{34}} x_{12} x_{26},\,
\underline{x_{34}} x_{15} x_{56},\,\,\,
\underline{x_{16}} x_{23} x_{35},\,
\underline{x_{16}} x_{24} x_{45},\,\,\,
\underline{x_{25}} x_{13} x_{14},\,
\underline{x_{25}} x_{36} x_{46}.
$$
The underlined variables are the non-edges
of the octahedron.

The abstract tight span $T$ has the following
geometric description. The polytope dual
to the truncated octahedron is gotten from
the $3$-cube by subdividing
each of the six facets
with a new vertex (thus creating
$6 \times 4 = 24$ edges) and then
erasing the $12$ edges of the cube.
Consider the six $4$-valent vertices
we just introduced. Each of them can be
replaced by two trivalent vertices
with a new edge in-between. If this replacement
is done in the most symmetric manner then
the result is a dodecahedron.
Finally, we glue a triangle and an edge
on each of the six new edges. The result is
Figure \ref{fig:dodec}.

\bigskip

\section{Prime Metrics and Minimal Subdivisions}

Koolen, Moulton and T\"onges \cite{KMT} classified all the 
\emph{prime metrics} for $n = 6$. These are the rays
in the metric fan $MF_6$. There are $14$ symmetry classes:

 1/5 Split: $\, d\,= \,(0, 0, 0, 0, 1, 0, 0, 0, 1, 0, 0, 1, 0, 1, 1) $

 2/4 Split: $\, d \, = (1, 1, 0, 1, 1, 0, 1, 0, 0, 1, 0, 0, 1, 1, 0) $

 3/3 Split: $\, d \, = (0, 0, 1, 1, 1, 0, 1, 1, 1, 1, 1, 1, 0, 0, 0) $

 Prime $P_1$ : $\, d \, = (1, 1, 1, 1, 2, 2, 2, 2, 1, 2, 2, 1, 2, 1, 1) $

 Prime $P_2$ : $\, d \, = (1, 1, 1, 2, 2, 2, 2, 1, 1, 2, 1, 1, 1, 1, 2) $

 Prime $P_3$ : $\, d \, = (1, 1, 1, 2, 2, 2, 2, 1, 1, 2, 1, 1, 3, 1, 2) $

 Prime $P_4$ : $\, d \, = (1, 1, 1, 1, 2, 1, 1, 1, 1, 1, 1, 1, 2, 1, 1) $

 Prime $P_5$ : $\, d \, = (1, 2, 2, 2, 4, 3, 3, 3, 3, 2, 2, 2, 4, 2, 2) $

 Prime $P_6$ : $\, d \, = (1, 1, 2, 3, 3, 2, 3, 2, 2, 3, 2, 2, 1, 1, 2) $

 Prime $P_7$ : $\, d \, = (1, 1, 1, 2, 2, 2, 2, 1, 1, 2, 1, 1, 2, 1, 2) $

 Prime $P_8$ : $\, d \, = (1, 2, 2, 4, 4, 3, 3, 3, 3, 4, 2, 2, 2, 2, 4) $

 Prime $P_9$ : $\, d \, = (1, 1, 1, 2, 3, 2, 2, 1, 2, 2, 1, 2, 3, 2, 1) $

 Prime $P_{10}$ : $\, d \, = (0, 1, 1, 1, 2, 1, 1, 1, 2, 2, 2, 1, 2, 1, 1) $

 Prime $P_{11}$ : $\, d \, = (0, 1, 1, 2, 2, 1, 1, 2, 2, 2, 1, 1, 1, 1, 2) $

\smallskip

\noindent Our computations provide an independent verification
of the correctness and completeness of the results in \cite{KMT}. Namely, we
computed the cone in the metric fan corresponding to each
of the $339$ metrics. For each cone we computed
(using POLYMAKE \cite{GJ}) the facets and the extreme rays
of the cone. And it turned out that the extreme rays
we found are precisely the $14$ types listed above.
All the facets and  extreme rays of 
the $339$ types of maximal cones in the 
metric fan $MF_6$ are posted at  
$$ \hbox{\tt bio.math.berkeley.edu/SixPointMetrics} $$
Our computations lead to the following result.

\begin{proposition} \label{fande}
If $f_i$ is the number of types of maximal cones in
the metric fan $MF_6$ with $i$ facets
and $e_j$ is the number with $j$ extreme rays then
\begin{eqnarray*}
\bigl( f_{15},f_{16},\ldots,f_{21}\bigr)   \quad = & \!\!\!
\bigl( 197, 42, 63, 18, 8, 10, 1  \bigr) \\
\bigl( e_{15}, e_{16}, \ldots, e_{24} \bigr)  \quad =  & \,\,\,
\bigl( 197, 60, 28, 19, 20, 2, 5, 2, 1, 5 \bigr).
\end{eqnarray*}
In particular, the metric fan $MF_6$ has $197$ types
of simplicial cones.
\end{proposition}

Proposition \ref{fande} shows that the largest number of facets 
of any cone is $21$, attained by only one type, and the largest 
number of extreme rays is $24$, attained by five types.  
The next two examples concern the extremal cases.

\begin{example}
Type 12 is the last metric listed in the table of Section 2.
Its cone in the metric fan is described by the following
 $21$ linear inequalities:
 \begin{eqnarray*}  &
 d_{12}  + d_{25} \geq d_{15} ,\,\,\,\,
  d_{13} + d_{36} \geq  d_{16},\,\,\,\,
  d_{45} + d_{46} \geq d_{56} ,\, \\ &
 d_{25} + d_{45} \geq d_{24} ,\,\,\,\,
 d_{36} + d_{46} \geq d_{34} ,\,\,\,\,
  d_{12} + d_{13} \geq d_{23} ,\, \\ &
 d_{26} + d_{34} \geq d_{24}  + d_{36} ,\,\,
  d_{15} + d_{26} \geq  d_{16} + d_{25},\, \,
 d_{14} + d_{23} \geq d_{13}  + d_{24} ,\,\\ &
  d_{14} + d_{56}  \geq d_{16} + d_{45},\, \,
  d_{16} + d_{35} \geq d_{15} + d_{36} ,\,\,
  d_{24}  + d_{35} \geq d_{25} + d_{34},\,\\ &
  d_{14} + d_{23} \geq d_{12} + d_{34} ,\,\,
  d_{14}  + d_{56} \geq d_{15} + d_{46} ,\, \,
  d_{26} + d_{34} \geq  d_{23} + d_{46} , \,\\ &
  d_{24} + d_{35}  \geq d_{23} + d_{45} ,\,\,
  d_{15} + d_{26} \geq  d_{12} + d_{56} , \,\,
  d_{16} + d_{35} \geq d_{13} + d_{56} ,\, \\ &
  d_{15} + d_{23}+ d_{34}  + d_{56} \geq  d_{16} + d_{24}  + 2 d_{35} ,\,\\ &
  d_{16} + d_{23} + d_{24}  + d_{56} \geq d_{15} + 2 d_{26} + d_{34} ,\,\\ &
  d_{15} + d_{16}  + d_{24} + d_{34} \geq  2 d_{14} + d_{23} + d_{56} .
\end{eqnarray*}
None of these $21$ inequalities is redundant. 
This cone has $19$ extreme rays: the six 1/5 splits,
six 2/4 splits,  three rays of type $P_6$, 
three rays of type $P_7$, and a unique
ray of type $P_2$, namely,
$(1, 1, 1, 2, 2, 2, 2, 1, 1, 2, 1, 1, 1, 1, 2)$. \qed
\end{example}

The six 1/6-splits are among the extreme rays
of every cone in the metric fan  $MF_6$. In the
next example, we only list the other extreme rays.

\begin{example} 
The five types of cones in $MF_6$ with $24$ extreme rays
are \# 7,  \# 26, \# 337, \# 338, and \# 339.
Only \# 7 corresponds to a two-dimensional tight span.
The following $18$ vectors are extreme rays of this cone:
$$
\begin{matrix}
         2/4 \,{\rm split}: (1 0 1 1 1 1 0 0 0 1 1 1 0 0 0) &
         2/4 \,{\rm split}: (1 0 0 0 1 1 1 1 0 0 0 1 0 1 1)  \\
         2/4 \,{\rm split}: (0 1 1 1 1 1 1 1 1 0 0 0 0 0 0)  &
         2/4 \,{\rm split}: (0 1 0 0 1 1 0 0 1 1 1 0 0 1 1)  \\
         P_2:     (1 1 1 2 2 2 2 1 1 2 1 1 1 1 2)  &
         P_6:     (1 2 1 3 3 3 2 2 2 3 1 1 2 2 2)  \\
         P_6:     (1 2 2 2 3 3 3 1 2 2 2 1 2 1 3)  &
         P_6:     (2 1 1 3 3 3 3 1 1 2 2 2 2 2 2)  \\
         P_6:     (2 1 2 2 3 3 2 2 1 3 1 2 2 1 3)  &
         P_7:     (1 1 1 2 2 2 2 1 1 2 1 1 2 1 2)  \\
         P_8:     (2 2 1 4 4 4 3 2 2 3 2 2 3 3 4)  &
         P_8:     (2 2 2 3 4 4 4 1 2 4 3 2 3 2 3)  \\

         P_8:     (2 2 2 3 4 4 4 3 2 4 1 2 3 2 3)  &
         P_8:     (2 2 3 4 4 4 3 2 2 3 2 2 3 1 4)  \\
         P_{10}:     (1 1 1 2 2 2 0 1 1 2 1 1 1 1 2) & 
         P_{10}:     (1 1 1 2 2 2 2 1 1 0 1 1 1 1 2)  \\
         P_{10}:     (1 1 1 0 2 2 2 1 1 2 1 1 1 1 2)  &
         P_{10}:     (1 1 1 2 2 2 2 1 1 2 1 1 1 1 0) 
\end{matrix}
$$
The other four types with $24$ extreme rays are three-dimensional,
and each of them has three-dimensional prime metrics
among its extreme rays. For instance, the following 
$18$ vectors are extreme rays of the cone \# 338:
$$
\begin{matrix}
   3/3 \,{\rm split}:    (1 1 0 0 1 0 1 1 0 1 1 0 0 1 1)  &
   3/3 \,{\rm split}:    (1 0 1 0 1 1 0 1 0 1 0 1 1 0 1)  \\
   3/3 \,{\rm split}:    (0 1 0 1 1 1 0 1 1 1 0 0 1 1 0)  &
   3/3 \,{\rm split}:    (0 0 1 1 1 0 1 1 1 1 1 1 0 0 0)  \\
   P_1:           (1 1 1 1 2 2 2 2 1 2 2 1 2 1 1)   &
   P_4:          (1 1 1 1 2 1 1 1 1 2 1 1 1 1 1)  \\
   P_5:           (3 2 2 2 4 3 3 3 1 4 2 2 2 2 2)  &
   P_5:           (2 2 2 1 4 2 2 3 2 4 3 2 3 2 3)  \\
   P_5:           (1 2 2 2 4 3 3 3 3 4 2 2 2 2 2)  &
   P_5:           (2 2 2 3 4 2 2 3 2 4 3 2 3 2 1) \\
   P_9:           (2 2 1 2 3 2 1 2 1 3 2 1 1 2 1) &
   P_9:           (1 2 1 1 3 1 2 2 2 3 1 1 2 2 2)  \\
   P_9:           (1 1 2 1 3 2 1 2 2 3 2 2 1 1 2)  & 
   P_9:           (2 1 2 2 3 1 2 2 1 3 1 2 2 1 1)  \\
   P_{11}:           (1 1 1 1 2 2 2 2 1 2 2 1 0 1 1)  &
   P_{11}:           (1 1 1 1 2 2 2 2 1 2 0 1 2 1 1)  \\
   P_{11}:           (1 1 1 1 2 0 2 2 1 2 2 1 2 1 1) &
   P_{11}:           (1 1 1 1 2 2 0 2 1 2 2 1 2 1 1) 
\end{matrix}
$$
The corresponding lists for all $339$ 
types appear on our website. \qed
\end{example}

It is instructive to draw the tight spans of the
eleven prime metrics. Four of them are
actually three-dimensional. For instance,
the metric $P_4$ has the structure of an octahedron.
Please compare  Figure 5 with  \cite[Figure 1]{KMT}.

For each of the $11$ prime metrics,
we list the $f$-vector of their tight span:
$$ \begin{matrix}
& & \# \,\,{\rm vertices} &
\# \,\,{\rm edges} &
\# \,\,{\rm polygons} &
\# \,\,\hbox{$3$-cells} \\
P_1: & & 6 & 9 & 4 & 0 \\
P_2: & & 7 & 15 & 9 & 0  \\
P_3: & & 6 & 10 & 6 & 1 \\
P_4: & & 10 & 20 & 12 & 1  \\
P_5: & & 11 & 20 & 11 & 1  \\
P_6: & & 7 & 11  & 5 & 0 \\
P_7: & & 11 & 20 & 10 & 0  \\
P_8: & & 11 & 19 &  9 & 0 \\
P_9: & &  7 & 12 &  7 & 1 \\
P_{10}: & & 5 &  7  & 3 & 0   \\
P_{11}: & & 5 &  7  & 3 & 0 \\
\end{matrix}
$$

The metrics $P_1,\ldots,P_9$ each 
have six \emph{trivial vertices} in their tight span.
The list of non-trivial vertices given in
\cite[Table 1]{KMT} is consistent with our
list above. The prime metrics $P_{10}$ and
$P_{11}$ are improper in the sense
that two points have distance zero.

\begin{figure}[htb]
\vskip -0.7cm
\epsfig{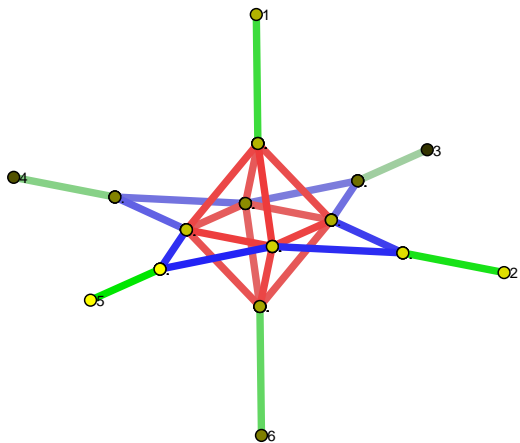}
\vskip -1.8cm
\caption{The tight span of the metric $P_4 + (1,1,\ldots,1)$}
\label{fig:prime4}
\end{figure}

\begin{remark}
The tight spans of the prime metrics $P_3,P_4,P_5,P_9$
have a $3$-dimensional cell, while the
tight spans of the other seven
are $2$-dimensional.
\end{remark}

The metric fan $MF_6$ defines an incidence
relation between the prime metrics and the generic metrics.
This leads to a finer invariant
for distinguishing among the $339$ types of generic types.
This invariant is the vector
$\,(S,P) = (s_2,s_3, \, p_1,p_2,\ldots,p_{11})\,$
where $p_i$ is the number of extreme rays
of type $P_i$ which lie on the corresponding cone.
This invariant resolves about half of the
clusters which had been left unresolved by
the earlier invariants.

\begin{example}
The three two-dimensional types \# 9, \# 10 and \# 11,
all have the same invariants $R= (2,8,5)$,
$ g = 2 $ and $ c=4$ in the list of Section 2.
They are distinguished by the new invariant
$\,(S,P) = (s_2,s_3, \, p_1,p_2,\ldots,p_{11})\,$

Type 9 has $\, (S,P) \, = \, 
(4, 0, 0, 1, 0, 0, 0, 2, 1, 2, 0, 2, 0) $.

Type 10 has $\, (S,P) \, = \, 
(5, 0, 0, 1, 0, 0, 0, 1, 2, 2, 0, 2, 0)$.

Type 11 has $\, (S,P) \, = \, 
(5, 0, 0, 1, 0, 0, 0, 2, 2, 2, 0, 2, 0)$.
\end{example}

\begin{example}
The last three three-dimensional types \# 337, \# 338 and \# 339
all have the same invariants $(f,R,B,S,C)$ 
in the big table of Section 3.
They are distinguished by the new invariant
$\,(S,P) = (s_2,s_3, \, p_1,p_2,\ldots,p_{11})\,$

Type 337 has $\, (S,P) \, = \, 
 (0, 4, 0, 0, 1, 1, 4, 0, 0, 0, 4, 0, 4) $.

Type 338 has $\, (S,P) \, = \, 
 (0, 4, 1, 0, 0, 1, 4, 0, 0, 0, 4, 0, 4)$.

Type 339 has $\, (S,P) \, = \, 
 (0, 4, 1, 0, 1, 0, 4, 0, 0, 0, 4, 0, 4)$.
\end{example}

We close this section with a remark aimed at experts in
polytope theory. The metric fan is the secondary
fan of the hypersimplex $\Delta(6,2)$, hence it
is the normal fan of the \emph{secondary polytope}
$\,\Sigma \bigl( \Delta(6,2) \bigr) $.
Following \cite{BFS}, the vertices of the secondary
polytope correspond to the regular triangulations
of $\Delta(6,2) $, and the facets of the secondary
polytope correspond to \emph{minimal regular subdivisions}
of $\Delta(6,2) $.
For instance, Figure 5 is the dual picture to
a regular subdivision of $\Delta(6,2)$ into $10$ 
five-dimensional polytopes.

The results described in this section provide the
vertex-facet incidence matrix of the secondary
polytope $\,\Sigma \bigl( \Delta(6,2) \bigr) $.
In particular, the classification result of Koolen, 
Moulton and T\"onges 
\cite{KMT} can be restated as follows.

\begin{corollary} Up to the action of the
symmetric group on $\{1,2,3,4,5,6\}$,
the secondary polytope $\,\Sigma \bigl( \Delta(6,2) \bigr) \,$
has precisely $14$~facets.
\end{corollary}

\section{Visualization of Tight Spans in POLYMAKE}

POLYMAKE is a software package
 developed by  Ewgenij Gawrilow and Michael Joswig for 
studying polytopes and polyhedra  \cite{GJ}.  
It allows us to define a polyhedron by a set of 
either inequalities or vertices and computes numerous 
properties of the polyhedron.  We implemented a client 
program to POLYMAKE for visualizing the tight span $T_d$ of a 
given metric $d$.  In short, our program does the 
following to produce the figures in this paper:
\begin{enumerate}
\item Compute  all faces of the polyhedron $P_d$ from the given metric $d$.
\item  Build the tight span $T_d$ by extracting the bounded faces of $P_d$.
\item Spring-embed the 2-skeleton of $T_d$ in  
$3$-space using POLYMAKE's spring embedder and display it using JAVAVIEW. 
\item Label the points corresponding to the finite set (the ``taxa'') on which the metric is defined. In our case, the set of labels is $\,\{1,2,3,4,5,6\}$.
\end{enumerate}
The program was developed by the second author with the assistance
of Michael Joswig and Julian Pfeifle.
The code can be downloaded at \ 
$$ \hbox{\tt www.math.berkeley.edu/$\sim$jyu/}.$$

The figures above only show the combinatorial properties of 
the tight span and, because of the projection and the spring embedder, 
the edge lengths seen here do not represent the actual edge 
lengths in the tight span.

The pictures produced by our software are
different from the output produced
by the software SPLITSTREE \cite{DHM}.
SPLITSTREE  is a program that can, among other things, 
compute and visualize the split-decompositions 
of metrics (see \cite{BD,DHM}).  
It decomposes an input metric into a sum of splits 
plus a \emph{split-prime} metric that cannot be 
further decomposed into splits.  SPLITSTREE outputs a planar 
graph representing the split-decomposable part of  the 
input metric, where sets of parallel edges represent splits.

When a metric is split-decomposable, then the tight span is a cubical
complex,  and the output from SPLITSTREE agrees with our visualization with a 
few edges removed to make it planar.
The edges output by SPLITSTREE are scaled to be metrically accurate.
By contrast, our implementation does not preserve the edge lengths.  
Among the $339$ generic metrics on six points, only one type
(namely, Type \# 66) is split-decomposable. For all other
$338$ types, the picture produced by  our program contains more
refined combinatorial information than the output of SPLITSTREE.

The split-decomposition theory of Bandelt and Dress
\cite{BD} has the following interpretation in terms
of the metric fan. The $31$ splits are among the
extreme rays of the metric fan. Consider the
\emph{induced subfan} of $MF_6$ on the
$31$ splits. A key result in \cite{BD}
states that this subcomplex is simplicial, i.e.,
all cones in this subfan are spanned by
linearly independent vectors. Now consider any
metric $d$, for instance, one of our $339$ generic metrics,
and let $C$ be the cone of 
the metric fan $MF_6$ containing $d$ in its
relative interior. The intersection of $C$ with  the induced
subfan of split-decomposable metrics is a simplicial face $F$ of $C$.
It follows that $d$ can be written uniquely
as a sum of a vector $d'$ in $F$ and a positive combination of
extreme rays of $C$ not in $F$. 
The output of SPLITSTREE is the edge graph
of the tight span  of $d'$.

\section{Acknowledgments}

This project grew out of the second
author's term project
for the Seminar on Mathematics of Phylogenetic Trees
which was organized by Lior Pachter and the first author
during the Fall Semester 2003 at UC Berkeley.
We are grateful to Michael Joswig,
Julian Pfeifle and J\"org Rambau for
helping us with our computations.
Josephine Yu was supported in part by
an NSF Graduate Research Fellowship. Bernd Sturmfels was supported by
a Hewlett Packard Visiting Research Professorship 2003/2004
at MSRI Berkeley and in part  by the NSF (DMS-0200729).

 \end{document}